\documentclass[12pt]{article}

\usepackage{amsxtra,amssymb,amsthm,amsmath,latexsym}
\usepackage{graphicx}
\usepackage{booktabs}
\usepackage{multirow}
\usepackage{afterpage}
\usepackage{setspace}
\usepackage{microtype}
\usepackage{calc,fourier}
\usepackage[T1]{fontenc}
\usepackage[hidelinks]{hyperref}
\usepackage{float}
\usepackage[small]{caption}
\allowdisplaybreaks

\newtheorem{thm}{Theorem}[]

\numberwithin{equation}{section}

\newcommand{\bee}{\begin{equation*}}
\newcommand{\eee}{\end{equation*}}
\newcommand{\be}{\begin{equation}}
\newcommand{\ee}{\end{equation}}
\newcommand{\ba}{\begin{align}}
\newcommand{\ea}{\end{align}}

\title{Simple Method for Evaluating Singular Integrals}
\author{N. T. Tran\footnote{Mailing address:  Mathematics Department, 138 Cardwell Hall, Manhattan, KS 66506} \\
\small Department of Mathematics\\
\small Kansas State University, Manhattan, KS 66506-2602, USA\\
\small \texttt{*nhantran@ksu.edu} 
}

\date{}
\begin{document}
\maketitle

\begin{abstract}
In this paper, we study the class of one dimensional singular integrals that converge in the sense of Cauchy principal value. In addition, we present a simple method for approximating such integrals. 
\end{abstract}



\section{Introduction} \label{sec1}
Many problems in engineering and science require evaluating singular integrals. For example, in electromagnetic and acoustic wave scattering, the boundary integral equations have singular kernels, see \cite{RMaterial, TEMImpedance, TFastScalar, TMaterial, TPerfect, TBook}. In fluid and solid mechanics, physicists and engineers face the same problem, see \cite{Chien, Karami}. Thus, the study of such integrals plays an important role in engineering and science. In this paper, we consider only one dimensional singular integrals that converge in the sense of Cauchy principal value.

One dimensional singular integrals are defined in the literature as follows
\begin{equation}
	\int_a^b \frac{u(t)}{(t-s)^p} dt, \quad s \in (a,b), \quad p>0,
\end{equation}
in which $u(t)$ is a continuous function. These integrals are classified by the order of singularity. If $p<1$, the integral is called weakly singular. If $p = 1$, the integral is strongly singular. If $p > 1$, the integral is called hyper-singular, see \cite{Neri}. In other words, an integral is called weakly singular if its value exists and continuous at the singularity. An integral is called strongly singular if both the integrand and integral are singular. An integral is called hyper-singular if the kernel has a higher-order singularity than the dimension of the integral. For strongly singular integrals, they are often defined in terms of Cauchy principal value, see \cite{Wazwaz}. For hyper singular integrals, they are often interpreted as Hadamard finite part integrals, see \cite{Hadamard}.

There are many special methods developed to treat singular integral problems since numerical integration routines often lead to inaccurate solutions. For example, to deal with the singularities in surface integral equations, the method of moments regularizes the singular integrals by sourcing them analytically for specific observation point \cite{Resende, Tzoulis}. Other methods include Gaussian quadrature method which has high-order of accuracy with a non-uniform mesh \cite{Hui, Tsamasphyros}, Newton-Cotes method which has low-order of accuracy with a uniform mesh \cite{Du, Li, Sun},  Guiggiani\textquotesingle s method which extracts the singular parts of the integrand and treat them analytically \cite{Guiggiani}, sigmoidal transformation which transforms the integrand to a periodic function \cite{Choi, Elliott}, and Duffy\textquotesingle s transformation which cancels the singularity of type $\frac{1}{t}$ \cite{Duffy}. Most of these methods can be characterized in three categories: singularity subtraction, analytical transformation, and special purpose quadrature.

In this paper, we present an alternative approach for approximating one dimensional singular integrals which converge in the sense of Cauchy principal value. In addition, a proof of this method is outlined in section \ref{sec2} to serve as a theoretical basis for the method. In section \ref{sec3}, the detailed implementation of our method is described for integrals over the standard interval [-1,1].

\section{Approximation of singular integrals} \label{sec2}
\begin{thm}
	Let $\int_D f(x)dx$, $D \subseteq [-1, 1]$, be a singular integral that has finite value in the sense of Cauchy principal value. Suppose $x_0$ is its singularity in $D$. Then, for any $\epsilon>0$, there exist $N>0$ and $a_j, 0\le j \le n$,  $n \ge N$, such that: for all $n \ge N$
	\be
	\left| \int_D f(x)dx - \sum_{j=0}^{n}a_j \int_{D} U_j(x)dx\right| < \epsilon,
	\ee
	where $U_j, 0\le j \le n$, are Chebyshev polynomials of second kind.
\end{thm}

\textbf{Proof}

Let $D_\delta:=D \setminus B_\delta(x_0)$, where $B_\delta(x_0):=(x_0-\delta, x_0+\delta)$ and $\delta>0$. Since $f$ is continuous in $D_\delta$, $f$ can be expressed as:
\be \label{eq1}
f(x) = \sum_{j=0}^{\infty} a_j U_j(x), \quad x \in D_\delta
\ee 
where $U_j$ are Chebyshev polynomials of second kind
\be \label{eq2}
U_j(x)=\frac{\sin((j+1)\cos^{-1}(x))}{\sin(\cos^{-1}(x))}, \quad  j \ge 0.
\ee
Therefore
\be \label{eq3}
\left| \int_{D_\delta} f(x)dx - \int_{D_\delta}\sum_{j=0}^{n}a_j U_j(x)dx\right| \to 0, \quad\text{ as } n\to \infty.
\ee
Since $\int_D f(x)dx$ has finite value in the sense of Cauchy principal value, one has
\be \label{eq4}
\left| \int_{D} f(x)dx - \int_{D_\delta} f(x)dx \right| \to 0 \quad\text{ as } \delta \to 0.
\ee
This means for any $\epsilon>0$, there exists $\delta_\epsilon>0$ such that: for all $0 < \delta < \delta_\epsilon$
\be \label{eq5}
\left| \int_{D} f(x)dx - \int_{D_\delta} f(x)dx \right| < \epsilon.
\ee
From \eqref{eq3}, for any $\epsilon$, there exists $N>0$ such that: for $n \ge N$
\be \label{eq6}
\left| \int_{D_\delta} f(x)dx - \int_{D_\delta}\sum_{j=0}^{n}a_j U_j(x)dx\right| < \epsilon.
\ee
Thus, from \eqref{eq5} and \eqref{eq6}
\be \label{eq7}
\left| \int_{D} f(x)dx - \int_{D_\delta}\sum_{j=0}^{n}a_j U_j(x)dx\right| < 2\epsilon,
\ee
for all $0 <\delta < \delta_\epsilon$, $n \ge N$. $\hfill\square$

\section{Methods for computing singular integrals} \label{sec3}
In this section, we present a method for evaluating the following singular integral which converges in the sense of Cauchy principal value
\be
S=\int_{-1}^{1} f(x)dx.
\ee
Without loss of generality, the singularity can be assumed to be at zero. For general case, one can always divide the interval of integration into many small intervals and treat them separately.

From section \ref{sec2}, we need to find the coefficient $a_i$ such that 
\be
S \simeq \sum_{i=0}^{n}a_i \int_{-1}^1 U_i(x)dx.
\ee

Since $U_i$ are Chebyshev polynomials of second kind, they admit some nice properties
\begin{align}
&\text{1. } \int_{-1}^1 U_i(x)dx=\frac{2\sin^2 \frac{(i+1)\pi}{2}}{i+1}, \\
&\text{2. } \int_{-1}^{1}U_j(x)U_i(x)\sqrt{1-x^2}dx=\frac{\pi}{2} \delta_{ij}.   	
\end{align}

Now consider the following integral
\begin{align}
\int_{-1}^{1} f(x)U_i(x)\sqrt{1-x^2} dx &\simeq \int_{-1}^{1} \sum_{j=0}^{n}a_j U_j(x)U_i(x)\sqrt{1-x^2} dx \\
&= \sum_{j=0}^{n} a_j\int_{-1}^{1} U_j(x)U_i(x)\sqrt{1-x^2} dx \\ 
&= \sum_{j=0}^{n} a_j \frac{\pi}{2} \delta_{ij} \\
&= \frac{\pi}{2} a_i, \quad 0 \le i \le n.
\end{align}
Thus, the coefficient $a_i$ can be computed by
\be
a_i = \frac{2}{\pi} \int_{-1}^{1} f(x)U_i(x)\sqrt{1-x^2} dx, \quad 0 \le i \le n,
\ee
and 
\begin{align}
S=\int_{-1}^{1} f(x)dx &\simeq \sum_{i=0}^{n}a_i \int_{-1}^1 U_i(x)dx \\
&= \sum_{i=0}^{n}a_i \frac{2\sin^2\frac{(i+1)\pi}{2}}{i+1}.
\end{align}

\section{Conclusions} \label{sec4}
In this paper, we present a method for approximating singular integrals which converge in the sense of Cauchy principal value. The proof of this method is outlined and the detailed implementation is also provided. One of the advantages of this method is that it is simple to implement. This method can serve as an alternative approach to other special methods in the literature.

\bibliographystyle{ieeetr}

\end{document}